\documentclass[12pt,a4]{article}
\usepackage{amsmath}
\usepackage{amsfonts}
\usepackage{euscript}
\def\scr{\EuScript}
\usepackage{amssymb}

\newcommand{\F}{\mathbb {F}}

\newtheorem{teorema}{Theorem.}[section]
\newtheorem{nota}[teorema]{Remark.}
\newtheorem{lema}[teorema]{Lemma.}

\newtheorem{proposicion}[teorema]{Proposition.}
\newtheorem{corolario}[teorema]{Corollary.}

\newenvironment{demostracion}{{\em Proof:\ }}{\hfill
 \rule{2mm}{2mm}  \par\vspace{2ex}}

\DeclareMathOperator{\height}{\rm ht}
\DeclareMathOperator{\Spec}{\rm Spec}
\newcommand{\tone}{t^{\frac{1}{p}}}
\newcommand{\uX}{\underline{X}}
\newcommand{\sP}{\scr{P}}
\newcommand{\cP}{\mathcal{P}}

\begin{document}

\date{March, 2001}
\title{Conservation of the noetherianity by perfect transcendental field extensions}

\author{M. Fern\'{a}ndez-Lebr\'{o}n and L. Narv\'{a}ez-Macarro\thanks{Both authors are partially supported by DGESIC, PB97-0723.}\\
Departamento de Algebra, Facultad de Matem\'aticas,\\ Universidad
de Sevilla, P.O. Box 1160, 41080 Sevilla, Spain.\\
lebron@algebra.us.es, narvaez@algebra.us.es}

\maketitle

\begin{abstract}
Let $k$ be a perfect field of characteristic $p>0$, $k(t)_{per}$
the perfect closure of $k(t)$ and $A$ a $k$-algebra. We
characterize whether the ring  $$A\otimes_k
k(t)_{per}=\bigcup_{m\geq 0}(A\otimes_k k(t^{\frac{1}{p^m}}))$$ is
noetherian or not. As a consequence, we prove that the ring
$A\otimes_k k(t)_{per}$ is noetherian when $A$ is the ring of
formal power series in $n$ indeterminates over $k$.
\end{abstract}


\section*{Introduction}

Motivated by the generalization of the results in \cite{nar_91}
(for  the case of a perfect base field $k$ of characteristic
$p>0$) in this paper we study the conservation of noetherianity by
the base field extension $k \rightarrow k(t)_{per}$, where
$k(t)_{per}$ is the perfect closure of $k(t)$. Since this
extension is not finitely generated, the conservation of
noetherianity is not clear {\em a priori} for $k$-algebras which
are not finitely generated.

Our main result states that $k(t)_{per} \otimes_k A$ is noetherian
if and only if $A$ is noetherian and for every prime ideal
${\mathfrak p}\subset A$ the field $\displaystyle \bigcap_{m\geq
0} Qt(A/{\mathfrak p})^{p^{m}}$ is algebraic over $k$ (see theorem
\ref{noet}). In particular, we are able to apply this result to
the case where $A$ is the ring of formal power series in $n$
indeterminates over $k$

We are indebted to J. M. Giral for giving us the proof of
proposition \ref{intersecc} and for other helpful comments.

\section{Preliminaries and notations} \label{sec-3}

All rings and algebras considered in this paper are assumed to be
commutative with unit element. If $B$ is a ring, we shall denote
by $\dim(B)$ its Krull dimension and by $\Omega (B)$ the set of
its maximal ideals. We shall use the letters $K,L,k$ to denote
fields and $\F _p$ to denote the finite field of $p$ elements, for
$p$ a prime number. If ${\mathfrak p} \in \Spec(B)$, we shall
denote by $\height({\mathfrak p})$ the height of ${\mathfrak p}$.
Remember that a ring $B$ is said to be {\em equicodimensional} if
all its maximal ideals have the same height. Also, $B$ is said to
be {\em biequicodimensional} if all its saturated chains of prime
ideals have the same length.

If $B$ is an integral domain, we shall denote by $Qt(B)$ its
quotient field.

 For any $\F_p$-algebra $B$, we denote
$B^{\sharp}={\displaystyle \bigcap_{m\geq 0} B^{p^m}}$.

We shall first study the contraction-extension process for prime
ideals relative to the ring extension $K[t]\subset K[\tone]$, $K$
being a field of characteristic $p>0$. Let us recall the following
well known result (cf. for example \cite{garling}, th. 10.8):

\begin{proposicion}\label{garling}
 Let $K$ be a field of characteristic $p>0$. Let $g(X)$ be a monic polynomial of
 $K[X]$. Then, the polynomial $f(X)=g(X^p)$ is irreducible in $K[X]$ if
and only if $g(X)$ is irreducible in $K[X]$ and not all its
coefficients are in $K^p$.
\end{proposicion}

From the above result, we deduce the following corollary.

\begin{corolario}\label{corgar}
Let $K$ be a field of characteristic $p>0$. Let $P$ be a non zero
prime ideal in $K[\tone]$ and let $F(t)\in K[t]$ be the monic
irreducible generator of the contraction $P^c = P\cap K[t]$. Then
the following conditions hold:
\begin{enumerate}
\item If $F(t)=a_0^p+a_1^pt+\cdots +t^d\in K^p[t]$, then
$P =(a_0+a_1\tone +\cdots +t^{\frac{d}{p}})$.
\item The equality $P=P^{ce}$ holds if
and only if $F(t)\not\in\ K^p[t]$.
\end{enumerate}
\end{corolario}

\begin{demostracion}

\begin{enumerate}

\item Consider the polynomial $G(\tau)=a_0+a_1\tau +\cdots
+\tau^d \in K[\tau] (\tau =\tone  )$ and the ring homomorphism
$\mu :K[\tau] \to K[t]$ defined by $$ \mu(\sum a_i \tau ^i )=\sum
a_i^p t^i.$$ From the identity $\mu(G) = F$ we deduce that
$G(\tau)$ is irreducible. Since $G(\tone)^p=F(t) \in P$, we deduce
that $G(\tone)\in P$ and then $P=(G(\tone))$.

\item The equality $P = P^{ce}$ means that $F(t) =
F(\tau^p)\in K[\tau]$ generates the ideal $P$, but that is
equivalent to saying that $F(\tau^p)$ is irreducible in $K[\tau]$.
To conclude, we apply proposition \ref{garling}.
\end{enumerate}
\end{demostracion}

For each $k$-algebra $A$, we define $A(t):=k(t)\otimes _k A$. We
also consider the field extension
$$k_{(\infty)}=\bigcup_{m\geq 1} k(t^{\frac{1}{p^m}}).$$ If $k$
is perfect, $k_{(\infty)}$ coincides with the perfect closure of
$k(t)$, $k(t)_{per}$.

 For the sake of brevity, we will write $t_m =
 t^{\frac{1}{p^m}}$. We also define
$$A_{(m)}:=A(t_m):=A\otimes_k k(t_m)=
  A(t)\otimes_{k(t)}k(t_m), \quad A_{[m]}:=A[t_m]$$
and $$A_{(\infty)}:=A\otimes_k k_{(\infty)}={\displaystyle
\bigcup_{m\geq 0} A_{(m)}},\quad A_{[\infty]}:={\displaystyle
\bigcup_{m\geq 0} A[t_m]}.$$

Each $A_{(m)}$ (resp. $A_{[m]}$) is a free  module over $A(t)$
(resp. over $A[t]$) of rank $p^m$ (because $(t_m)^{p^m}-t=0$).

For each prime ideal $P$
 of $A_{(\infty)}$ we denote $P_{[\infty]}:=P\cap A_{[\infty]}$,
 $P_{[m]}:=P\cap A_{[m]} \in \Spec(A_{[m]})$ and
 $P_{(m)}:=P\cap A_{(m)} \in \Spec(A_{(m)})$.

 In a similar way, if $Q$ is a prime ideal of $A_{[\infty]}$
we denote $Q_{[m]}:=Q\cap A_{[m]} \in \Spec(A_{[m]})$.

 We have:
\begin{itemize}
 \item
${\displaystyle P=\bigcup_{m\geq 0}P_{(m)}}$, ${\displaystyle
P_{[\infty]}=\bigcup_{m\geq 0} P_{[m]}}$, (resp. ${\displaystyle
Q=\bigcup_{m\geq 0}Q_{[m]}}$).
 \item $P_{(n)}\cap A_{(m)}=P_{(m)}$ and $P_{[n]}\cap A_{[m]}=P_{[m]}$ for all  $n\geq
m$ (resp. $Q_{[n]}\cap A_{[m]}=Q_{[m]}$ for all  $n\geq m$).
\end{itemize}

The following properties are straightforward:

\begin{enumerate}
\item The $k$-algebras $A_{[m]}$ (respectively $A_{(m)}$) are
 isomorphic to each other.
\item If $S_m = k[t_m] -\{ 0\}$, then $A_{(m)}=
 S_m^{-1} A_{[m]}$.
\item Since  $(S_m)^{p^m}\subset S_0\subset S_m $, we have
$A_{(m)}=S_0^{-1} A_{[m]}$ for $m\geq 0$. Consequently
$A_{(\infty)}=S_0^{-1} A_{[\infty]}$.
\item If $A$ is a domain (integrally closed), then $A_{[m]}$ and $A_{(m)}$ are
domains (integrally closed) for all $m\geq 0$ or $m =\infty$.
\item If $A$ is a noetherian k-algebra, then $A_{[m]}$ and
$A_{(m)}$ are noe\-the\-rian rings, for every $m\geq 0$.
\item If $A=k[\uX]=k[X_1,\dots ,X_n]$, then
 $A_{[\infty]}$ is not noetherian (the ideal generated by the
 $t_m$, $m\geq 0$, is not finitely generated).
\item If $I\subset A$ is an ideal, then $(A/I)_{(\infty)} =
 A_{(\infty)}/A_{(\infty)}I$.
\item If $T\subset A$ is a multiplicative subset, then
$(T^{-1}A)_{(\infty)}= T^{-1}A_{(\infty)}$.
\item If $A=k[\uX]$, then $A_{(\infty)}=k_{(\infty)}[\uX]$, hence
$A_{(\infty)}$ is noetherian. Moreover, $A_{(\infty)}$ is
noetherian for every finitely generated $k$-algebra $A$.
\end{enumerate}

The main goal of this paper is to characterize whether the ring
$A_{(\infty)}$ is noetherian (see th. \ref{noet} and corollary
\ref{noet-series}).

\begin{proposicion} \label{buena} With the above
notations, the following properties hold:
  \begin{enumerate}
   \item The extensions $A_{[m-1]}\subset A_{[m]}$ and $A_{(m-1)}\subset A_{(m)}$
   are finite free, and therefore integral and faithfully flat.
   \item The corresponding extensions to their quotient fields are purely
    inse\-pa\-ra\-ble.
  \end{enumerate}
\end{proposicion}
\begin{demostracion}
Straightforward.
\end{demostracion}

\begin{corolario}\label{entera}
$A_{[\infty]}$ (resp. $A_{(\infty)}$) is integral and faithfully
flat over each $A_{[m]}$ (resp. over each $A_{(m)}$).
\end{corolario}

From the properties above, we obtain the following lemmas:

\begin{lema}\label{prim}
Let $P'\subseteq P$ be prime ideals of $A_{(\infty)}$ (resp. of
$A_{[\infty]}$). The following conditions are equivalent:
\begin{enumerate}
 \item[(a)] $P'\subsetneq P$
 \item[(b)] There exists an $m\geq 0$ such that $P'_{(m)}\subsetneq P_{(m)}$
 (resp. $P'_{[m]}\subsetneq P_{[m]}$).
 \item[(c)] For every $m\geq 0$, $P'_{(m)}\subsetneq P_{(m)}$
(resp. $P'_{[m]}\subsetneq P_{[m]}$).
\end{enumerate}
\end{lema}

\begin{lema}\label{max}
Let $P$ prime ideal of $A_{(\infty)}$ (resp. of $A_{[\infty]}$).
The following conditions are equivalent:
\begin{enumerate}
 \item[(a)] $P$ is maximal.
 \item[(b)] $P_{(m)}$ (resp. $P_{[m]}$) is maximal for
 some $m\geq 0$.
 \item[(c)] $P_{(m)}$ (resp. $P_{[m]}$) is maximal for
 every $m\geq 0$.
\end{enumerate}
\end{lema}

\begin{corolario}  \label{ext-1}
With the notations above, for every prime ideal $P$ of
$A_{(\infty)}$ we have
$\height(P)=\height(P_{(m)})=\height(P_{[m]})$ for all $m\geq 0$.
Moreover, $\dim (A_{(\infty)}) = \dim(A_{(m)})$.
\end{corolario}

\begin{demostracion}
Since flat ring extensions satisfy the ``going down" property,
corollary \ref{entera} implies that $\height (P\cap A_{(m)})\leq
\height (P)$. By corollary \ref{entera} again, $A_{(\infty)}$ is
integral over $A_{(m)}$, then $\height (P) \leq \height (P\cap
A_{(m)})$.

 The equality $\height(P_{(m)})=\height(P_{[m]})$ comes
from the fact that $A_{(m)}$ is a localization of $A_{[m]}$.

The last relation is a standard consequence of the ``going up"
property.
\end{demostracion}

\begin{nota} \label{nota-ext-1} Corollary \ref{ext-1} remains
true if we replace $A_{(m)}\subset
A_{(\infty)}$ by $A_{[m]}\subset A_{[\infty]}$.
\end{nota}

\begin{corolario}  \label{ext-2}
With the notations above, for every  $Q\in \Spec(A_{(m)})$ there
is a unique $\widetilde{Q}\in\Spec(A_{(m+1)})$ such that
$\widetilde{Q}^c = Q$. Moreover, the ideal $\widetilde{Q}$ is
given by $\widetilde{Q} = \{y\in A_{(m+1)}\ |\ y^p\in Q\}$.
\end{corolario}

\begin{demostracion} This is an easy consequence of the fact that
$(A_{(m+1)})^p \subset A_{(m)}$.
\end{demostracion}

\begin{corolario}\label{cor}
Let us assume that $A$ is noetherian and for every maximal ideal
${\mathfrak m}$ of $A$, the residue field $A/{\mathfrak m}$ is
algebraic over $k$. Then for every $m\geq 0$ we have:
\begin{enumerate}
\item $\dim(A_{[\infty]})=\dim(A_{[m]})=\dim(A[t])=n+1$.
\item $\dim(A_{(\infty)})=\dim(A_{(m)})=\dim(A(t))=n$.
\end{enumerate}
\end{corolario}

\begin{demostracion}
The first relation comes from remark \ref{nota-ext-1} and the
noetherianity hypothesis.

The second relation comes from corollary \ref{ext-1} and
proposition (1.4) of \cite{nar_91}.
\end{demostracion}

The following result is a consequence of theorem (1.6) of
\cite{nar_91}, lemma \ref{max} and corollary \ref{cor}.

\begin{corolario} \label{cor-2}
  Let $A$ be a noetherian, biequidimensional, universally
ca\-te\-na\-rian $k$-algebra of Krull dimension $n$, and that for
any maximal ideal ${\mathfrak m}$ of $A$, the residue field
$A/{\mathfrak m}$ is algebraic over $k$. Then every maximal ideal
of $A_{(\infty)}$ has height $n$.
\end{corolario}

\section{The biggest perfect subfield of a formal functions field}
\label{sec-4}

Throughout this section, $k$ will be a perfect field of
characteristic $p>0$, $A=k[[\uX]]$, ${\mathfrak p}\subset A$ a
prime ideal, $R=A/{\mathfrak p}$ and $K=Qt(R)$.

The aim of this section is to prove that the biggest perfect
subfield of $K$, $K^{\sharp}= {\displaystyle \bigcap_{e\geq 0}
K^{p^e}}$, is an algebraic extension of the field of constants,
$k$. This result is proved in prop. \ref{intersecc} and it is one
of the ingredients in the proof of corollary \ref{noet-series}.

\begin{proposicion}\label{int1}
Under the above hypothesis, it follows that $k=R^{\sharp}$.
\end{proposicion}
\begin{demostracion}
Let ${\mathfrak m}$ be the maximal ideal of $R$. It suffices to
prove that $R^{\sharp}\subseteq k$. If $f\in R^{\sharp}$, then for
every $e>0$ there exists an $f_e\in R$ such that $f=f_e^{p^e}$.
\begin{itemize}
\item Suppose at first that $f$ is not a unit, then $f_e$ is not a unit for any
  $e>0$, and $f_e\in {\mathfrak m}$ for every $e>0$. Thus, $f\in {\mathfrak m}^{p^e}$
  for every $e>0$ and by Krull's intersection theorem,
  $$f\in \bigcap_{e\geq 0} {\mathfrak m}^{p^e}=\bigcap_{r\geq 0} {\mathfrak
  m}^r=(0).$$
\item If $f$ is unit, then $f=f_0+\widetilde{f}$, with $f_0\in
  k\subset R^{\sharp}$ and $\widetilde{f}\in R^{\sharp}$ and $f_0$ is
  unit. By the above case $\widetilde{f}=0$, hence $f\in k$.
\end{itemize}
\end{demostracion}

\begin{proposicion}\label{p0}

 If ${\mathfrak p}=(0)$, that is $R=k[[\uX]]$, $K=k((\uX))$, then $k=K^{\sharp}$.
\end{proposicion}
\begin{demostracion} It is a consequence of prop. \ref{int1} and
the fact that $R$ is a unique factorization domain.
\end{demostracion}

 In order to treat the general case, let us look at some
general lemmas.

\begin{lema}\label{alg}
(cf. \cite{bou_a_4_7} Chap. 5, $\S$ 15, ex. 8) If $L$ is a
separable algebraic extension of a field $K$ of characteristic
$p>0$, then $L^{\sharp}$ is an algebraic extension of
$K^{\sharp}$.
\end{lema}
\begin{demostracion}
 If $x\in L^{\sharp}$, then $x=y_e^{p^e}$ with $y_e\in L$ for all $e\geq
0$. Since $y_e$ is separable over $K$, $K(y_e)=K(y_e^{p^e})=K(x)$,
it follows that $y_e=x^{p^{-e}}\in K(x)$ and then
 $x\in K^{p^e}(x^{p^e})$. Therefore
$$[K^{p^e}(x):K^{p^e}]=[K^{p^e}(x^{p^e}):K^{p^e}]=[K(x):K].$$ Thus
$x$ satisfies the same minimal polynomial over $K^{p^e}$ and over
$K$ for all $e\geq 0$, and the coefficients of this minimal
polynomial must be in $K^{\sharp}$. So $x$ is algebraic over
$K^{\sharp}$.
\end{demostracion}

\begin{lema}\label{per}
Every algebraic extension of a perfect field is perfect.
\end{lema}
\begin{demostracion}
This is obvious because this is true for the finite algebraic
extensions.
\end{demostracion}

\begin{lema}\label{ej8}
Let $C$ be a subring of a domain $D$ and let $\overline{C}$ be the
integral closure of $C$ in $D$. If $f(X),g(X)$ are monic
polynomials in $D[X]$ such that $f(X)g(X)\in \overline{C}[X]$,
then $f(X),g(X)\in \overline{C}[X]$.
\end{lema}
\begin{demostracion}
We consider a field $L$ containig $D$ such that the polynomials
$f(X),f(X)$ are a product of linear factors:
$f(X)=\prod(x-\alpha_i)$, $g(X)=\prod(x-\beta_j)$,
$\alpha_i,\beta_j\in L$. Each $\alpha_i$ and $\beta_j$ are roots
of $f(X)g(X)$, hence they are integral over $\overline{C}$. Thus
the coefficients of $f(X)$ and $g(X)$ are integral over
$\overline{C}$ and therefore they are in $\overline{C}$.
\end{demostracion}

\begin{proposicion}\label{intersecc} Let $k$ be a perfect field of
characteristic $p>0$, $A=k[[\uX]]= k[[X_1,\dots,X_n]]$,
${\mathfrak p}\subset A$ a prime ideal, $R= A/{\mathfrak p}$  and
$K=Qt(R)$. Then $K^{\sharp}$ is an algebraic extension of $k$.
\end{proposicion}
\begin{demostracion}\footnote{Due to J. M. Giral.} Let
$r=\dim(A/{\mathfrak p})\leq n$. By the normalization lemma for
power series rings (cf. \cite{abhy_64}, 24.5 and
23.7)\footnote{The proof of the normalization lemma for power
series rings in \cite{abhy_64} uses generic linear changes of
coordinates and needs the field $k$ to be infinite. This proof can
be adapted for an arbitrary perfect coefficient field (infinite or
not) by using non linear changes of the form
$Y_i=X_i+F_i(X_{i+1}^p,\dots ,X_n^p)$, where the $F_i$ are
polynomials with coefficients in $\F_p$.}, there is a new system
of formal coordinates $Y_1,\dots,Y_n$ of $A$, such that
\begin{itemize}
 \item ${\mathfrak p}\cap k[[Y_1,\dots,Y_r]]=\{0\}$,
 \item $k[[Y_1,\dots ,Y_r]]\hookrightarrow {\displaystyle
  \frac{A}{{\mathfrak p}}}=R$ is a finite extension, and
 \item $k((Y_1,\dots,Y_r))\hookrightarrow K $ is a separable
 finite extension.
\end{itemize}
\noindent The proposition is then a consequence of proposition
\ref{p0} and lemma \ref{alg}\footnote{In particular, if $k$ is
algebraically closed, we would have $K^{\sharp}=k$.}.
\end{demostracion}

\begin{nota} Actually, under the hypothesis of proposition \ref{intersecc},
J.M. Giral and the authors have proved that the following stronger
properties hold: \begin{enumerate} \item[(1)] If $R$ is integrally
closed in $K$, then $K^{\sharp}=k$. \item[(2)] In the general
case, $K^{\sharp}$ is a \underline{finite} extension of
$k$.\end{enumerate}
\end{nota}

\section{Noetherianity of $A\otimes_k k(t)_{per}$}

Throughout this section, $k$ will be a perfect field of
characteristic $p>0$, keeping the notations of section
\ref{sec-3}.

\begin{proposicion}\label{p1}
Let $K$ be a field extension of $k$ and suppose that $K^\sharp$ is
algebraic over $k$. For every prime ideal
$\sP\in\Spec(K_{[\infty]})$ such that $\sP\cap k[t] = 0$ there
exists an $m_0\geq 0$ such that $\sP_{[m]}$ is the extended ideal
of $\sP_{[m_0]}$ for all $m\geq m_0$.
\end{proposicion}

\begin{demostracion} The extension $k[t]\subset K^{\sharp}[t]$ is
integral and then $\sP\cap K^{\sharp}[t]=0$.

We can suppose $\sP\neq (0)$. From Remark \ref{nota-ext-1}, we
have $\height(\sP_{[i]})=\height(\sP)= 1$ for every $i\geq 0$. Let
$F_i(t_i) \in K[t_i]$ be the monic irreducible generator of
$\sP_{[i]}$. From \ref{corgar}, for each $i\geq 0$ there are two
possibilities:
\begin{enumerate}
\item[(1)] $F_i\in K^p[t_i]$, then  $F_{i+1}(t_{i+1}) = F_i(t_i)^{1/p}$.
\item[(2)] $F_i\notin K^p[t_i]$, then $\sP_{[i+1]}=(\sP_{[i]})^{e}$
and $F_{i+1}(t_{i+1}) = F_i(t_i) = F_i(t_{i+1}^p)$.
\end{enumerate}

Since $\sP \cap K^{\sharp}[t]=(0)$, $F_0(t_0)\notin
({\displaystyle \bigcap_{m\geq 0}K^{p^m})[t_0]=\bigcap_{m\geq 0}
K^{p^m}[t_0]}$ and there exists an $m_0\geq 0$ such that
$F_0(t_0)\in K^{p^{m_0}}[t_0]$ and $F_0(t_0)\notin
K^{p^{m_0+1}}[t_0]$.

From (1) we have $F_i(t_i) = F_0(t_0)^{1/p^i} \in
K^{p^{m_0-i}}[t_i]$ for $i=0,\dots,m_0-1$ and
$F_{m_0}(t_{m_0})\notin K^p[t_{m_0}]$. Hence, applying (2)
repeatedly we find $F_{j+m_0}(t_{j+m_0}) = F_{m_0}(t_{m_0}) =
F_{m_0}(t_{j+m_0}^{p^j})$  and $\sP_{[j+m_0]}$ is the extended
ideal of $\sP_{[m_0]}$ for all $j\geq 1$.
\end{demostracion}

\begin{corolario}\label{cp1} Under the same hypothesis of
proposition \ref{p1}, $\sP$ is the extended ideal of some
$\sP_{m_0}$.
\end{corolario}
\begin{demostracion}
This is a consequence of prop. \ref{p1} and the equality
$\sP={\displaystyle \bigcup_{m\geq 0} \sP_m}$.
\end{demostracion}

Let $B$ be a free algebra over a ring $A$ and $S\subset A$ a
multiplicative subset. We denote by $I \mapsto I^E, J\mapsto J^C$
(resp. $I\mapsto I^e, J\mapsto J^c$) the extension-contraction
process between the rings $A$ or $S^{-1}A$ (resp. $A$ or $B$) and
the rings $B$ or $S^{-1}B$ (resp. $S^{-1}A$ or $S^{-1}B$).

\begin{proposicion}\label{p2}
With the notations above, let $\cP_1$ be a prime ideal in $B$ such
that $\cP_1\cap S =\emptyset$. Let $\cP_0 = \cP_1^C$, $\sP_1 =
\cP_1^e$ and $\sP_0 = \sP_1^C$. If $\sP_1=\sP_0^{E}$, then
$\cP_1=\cP_0^{E}$.
\end{proposicion}

\begin{demostracion} Let $\{ e_i\}$ be a $A$--basis of $B$. Since
$\cP_1\cap S = \emptyset$, it is clear that $\sP_1^c = \cP_1$,
$\sP_0^c = \cP_0$ and $\sP_0 = \cP_0^e$. If $\sP_1=\sP_0^{E}$, we
have
 $$\cP_1 = \cP_1^{ec} = \sP_1^c = (\sP_0^{E})^c = (\cP_0^{eE})^c
 =(\cP_0^{Ee})^c = (\cP_0^{E})^{ec} = \sum_{s\in S} (\cP_0^{E} :
 s)_{B}\supset \cP_0^{E}.$$ To prove the other inclusion,
 take an $s\in S$ and let $f=\sum a_ie_i$ be an element of
 $(\cP_0^{E} :s)_{B}$ with $a_i\in A$. Then, ${\displaystyle sf=\sum
 (sa_i)e_i\in \cP_0^{E}}$ and from the equality
 $\cP_0^{E}=\{\sum b_ie_i \ | \ b_i\in \cP_0\}$ we deduce that
 $sa_i\in \cP_0$ and $a_i\in (\cP_0^{E} :s)_{A}=\cP_0$.
 Therefore $f\in \cP_0^{E}$.

\end{demostracion}

\begin{proposicion}\label{p3}
Let $R$ be an  integral $k$-algebra, $K=Qt(R)$, and suppose that
$K^\sharp$ is algebraic over $k$. Then any prime ideal ${\cP}\in
\Spec(R_{[\infty]})$ with ${\cP}\cap k[t]=0$ and ${\cP} \cap R=0 $
is the extended ideal of some $\cP_{[m_0]}$, $m_0\geq 0$.
\end{proposicion}

\begin{demostracion}
Let us write $T=R-\{ 0\}$. We have $K= T^{-1}R$ and $K_{[m]} =
T^{-1} R_{[m]}$ for all $m\geq 0$ or $m=\infty$. We define $\sP =
T^{-1} \cP$. We easily deduce that $\sP_{[m]} = T^{-1} \cP_{[m]}$
for all $m\geq 0$.

From proposition \ref{p1}, there exists an $m_0\geq 0$ such that
$\sP_{[m]}$ is the extended ideal of ${\sP}_{[m_0]}$ for every
$m\geq m_0$. Then, proposition \ref{p2} tells us that $\cP_{[m]}$
is the extended ideal of ${\cP}_{[m_0]}$ for every $m\geq m_0$, so
$\cP = \bigcup \cP_{[m]}$ is the extended ideal of
${\cP}_{[m_0]}$.
\end{demostracion}

\begin{proposicion} \label{not-noet}
Let $K$ be a field extension of $k$ and suppose that $K^\sharp$ is
not algebraic over $k$. Then $K_{(\infty)}$ is not noetherian.
\end{proposicion}

\begin{demostracion} Let $s\in K^\sharp$ be a transcendental
element over $k$.

For each $m\geq 0$, let $s_m = s^{\frac{1}{p^m}}\in K$ and
$\alpha_m = t_m- s_m$. Let $P$ be the ideal in $K_{(\infty)}$
generated by the $\alpha_m, m\geq 0$. We have $\alpha_m =
\alpha_{m+1}^p$ and $P_{(m)} = K_{(m)}\alpha_m$ for all $m\geq 0$.

Suppose that $P$ is finitely generated. Then, there exists an
$m_0\geq 0$ such that $P = K_{(\infty)}\alpha_{m_0}$. By faithful
flatness, we deduce that $\alpha_{m_0+1}\in
K_{(m_0+1)}\alpha_{m_0}$. Let us write $\tau = t_{m_0+1}, \sigma
=s_{m_0+1}$. Then, $\alpha_{m_0+1}=\tau-\sigma $ and there exist
$\psi(\tau)\in K[\tau]= K_{[m_0+1]}$, $\varphi(\tau)\in
k[\tau]\setminus\{0\}$ such that $$ \varphi(\tau) (\tau-\sigma)=
\psi(\tau) (\tau-\sigma)^p.$$Simplifying and making $\tau=\sigma$
we obtain $$\varphi(\sigma) = \psi(\sigma) (\sigma-\sigma)^{p-1} =
0$$contradicting the fact that $s$ is transcendental over $k$.

We conclude that $P$ is not finitely generated and $K_{(\infty)}$
is not noetherian.
\end{demostracion}

\begin{teorema}\label{noet}
Let $k$ be a perfect field of characteristic $p>0$ and let $A$ be
a $k$-algebra. The following properties are equivalent:
\begin{enumerate}
\item[(a)] The ring $A$ is noetherian and for any ${\mathfrak p}\in\Spec(A)$,
the field $Qt(A/{\mathfrak p})^{\sharp}$ is algebraic over $k$.
\item[(b)] The ring $A_{(\infty)}$ is noetherian.
\end{enumerate}
\end{teorema}

\begin{demostracion} Let first prove (a) $\Rightarrow$ (b).
By Cohen's theorem (cf. \cite{naga_75}, (3.4)), it is enough to
prove that any $P\in \Spec(A_{(\infty)})-\{ (0)\}$ is finitely
generated.

From corollaries \ref{ext-1} and \ref{cor}, we have
$$\height(P_{[m]})=
\height(P_{(m)})=\height(P_{[\infty]})=\height(P)=r\leq n .$$
Consider the prime ideal of $A$:
 $${\mathfrak p}:= A\cap P = A\cap P_{[\infty]} = A\cap P_{[m]} = A \cap
 P_{(m)}.$$
There are two possibilities (cf. \cite{ega_iv_2}, prop. (5.5.3)):
\begin{enumerate}
 \item[(i)] $\height({\mathfrak p})=r=\height(P_{[m]})$ and $P_{[m]}={\mathfrak
p}[t_m]$, for every $m\geq 0$.
 \item[(ii)] $\height({\mathfrak p})=r-1=\height(P_{[m]})-1$,
${\mathfrak p}[t_m]\varsubsetneq P_{[m]}$ and $A/{\mathfrak p}
\varsubsetneq A[t_m]/P_{[m]}$ is algebraic generated by $t_m$ mod
$P_{[m]}$, for every $m\geq 0$.
\end{enumerate}

In case (i), $P_{[\infty]}$ and $P$ are the extended ideals of
$\mathfrak p$ and they are finitely generated.

Suppose we are in case (ii). We denote $R=A/{\mathfrak p}$,
$K=Qt(R)$.

Then: $$R_{[m]} = A_{[m]}/{\mathfrak p}[t_m], \quad R_{[\infty]} =
A_{[\infty]}/A_{[\infty]}{\mathfrak p}=
A_{[\infty]}/\bigcup_{m\geq 0}{\mathfrak p}[t_m].$$ Define
${\displaystyle \cP := R_{[\infty]} P_{[\infty]} =
P_{[\infty]}/\bigcup_{m\geq 0}{\mathfrak
p}[t_m]\in\Spec(R_{[\infty]})}$. We have $\cP_{[m]} =\cP\cap
R_{[m]}= P_{[m]}/{\mathfrak p}[t_m]$, $\cP\cap R = \cP\cap k[t] =
0$ and $$\height(\cP_{[m]})=\height\left( P_{[m]}/{\mathfrak
p}[t_m]\right) =1, \quad \height (\cP)= \height\left(
P_{[\infty]}/\bigcup_{m\geq 0}{\mathfrak p}[t_m]\right) =1.$$ We
conclude by applying proposition \ref{p3}: there exists an
$m_0\geq 0$ such that $\cP$ is the extended ideal of
$\cP_{[m_0]}$. Then, $P_{[\infty]}$ is the extended ideal of
$P_{[m_0]}$ and $P=A_{(\infty)}P_{[\infty]}=A_{(\infty)}P_{[m_0]}$
is finitely generated.

Let us prove now (b) $\Rightarrow$ (a). Since $A_{(\infty)}$ is
faithfully flat over $A$, we deduce that $A$ is noetherian.

Let ${\mathfrak p}\in \Spec(A)$ and let $R=A/{\mathfrak p}$,
$K=Qt(R)$. Noetherianity of $A_{(\infty)}$ implies, first,
noetherianity of $R_{(\infty)}$, and second, noetherianity of
$K_{(\infty)}$. To conclude we apply proposition \ref{not-noet}.
\end{demostracion}

\begin{corolario} \label{cor-2bis} Let $k$ be a perfect field of characteristic $p>0$ and let $A$ be
a noetherian $k$-algebra. The following properties are equivalent:
\begin{enumerate}
\item[(a)] The ring $A_{(\infty)}$ is noetherian.
\item[(b)] The ring $\left(A_{\mathfrak m}\right)_{(\infty)}$ is
noetherian for any maximal ideal ${\mathfrak m}\in \Omega(A)$.
\end{enumerate}
\end{corolario}

\begin{demostracion} For (a) $\Rightarrow$ (b) we use the fact
that $\left(A_{\mathfrak m}\right)_{(\infty)} = A_{\mathfrak
m}\otimes_A A_{(\infty)}$.

For (b) $\Rightarrow$ (a), let ${\mathfrak p}\subset A$ be a prime
ideal and let $\mathfrak m$ be a maximal ideal containing
$\mathfrak p$. From hypothesis (b) , the ring $\left(A_{\mathfrak
m}\right)_{(\infty)}$ is noetherian. Then, from theorem \ref{noet}
we deduce that the field $Qt(A/{\mathfrak
p})^\sharp=Qt\left(A_{\mathfrak m}/ A_{\mathfrak m} {\mathfrak
p}\right)^\sharp $ is algebraic over $k$. From theorem \ref{noet}
again we obtain (a).
\end{demostracion}

\begin{corolario}\label{noet-series}
Let $k$ be a perfect field of characteristic $p>0$, $k'$ an
algebraic extension of $k$ and $A=k'[[X_1,\dots,X_n]]$. Then, the
ring $A_{(\infty)} = k(t)_{per}\otimes_k A$ is noetherian.
\end{corolario}

\begin{demostracion} It is a consequence of lemma \ref{per},
proposition \ref{intersecc} and  theorem \ref{noet}.
\end{demostracion}

\begin{corolario} \label{cor-3}
Let $k$ be a perfect field of characteristic $p>0$. If
$(B,{\mathfrak m})$ is a local noetherian $k$-algebra such that
$B/{\mathfrak m}$ is algebraic over $k$, then
$B_{(\infty)}=k(t)_{per}\otimes_k B$ is noetherian. In particular,
the field $Qt(B/{\mathfrak p})^\sharp$ is algebraic over $k$ for
every prime ideal ${\mathfrak p}\subset B$.
\end{corolario}
\begin{demostracion} Let $k' = B/{\mathfrak m}$.
By Cohen structure theorem (cf. \cite{ega_iv_2}, Chap. 0, Th.
(19.8.8)), the completion $\widehat{B}$ of $B$ is a quotient of a
power-series ring $A$ with coefficients in $k'$. Since
$\widehat{B}_{(\infty)}$ is also a quotient of $A_{(\infty)}$, we
deduce from corollary \ref{noet-series} that $\hat{B}_{(\infty)}$
is noetherian. Since $\widehat{B}$ is faithfully flat over $B$,
the ring $\widehat{B}_{(\infty)}$ is also faithfully flat over
$B_{(\infty)}$. So, $B_{\infty}$ is noetherian.

The last assertion is a consequence of theorem \ref{noet}.
\end{demostracion}

\begin{corolario} Let $k$ be a perfect field of characteristic
$p>0$. For any noetherian $k$-algebra $A$ such that the residue
field $A/{\mathfrak m}$ of every maximal ideal ${\mathfrak
m}\in\Omega(A)$ is algebraic over $k$, the ring $A_{(\infty)}$ is
noetherian. Furthermore, if $A$ is regular and equicodimensional
then $A_{(\infty)}$ is also regular and equicodimensional of the
same dimension as $A$.
\end{corolario}

\begin{demostracion} The first part is a consequence of corollaries
\ref{cor-2bis} and \ref{cor-3}. For the last part, we use
corollary \ref{cor-2}, the fact that all $A_{(m)}, m\geq 0$ are
regular and of the same (global homological = Krull) dimension
(\cite{nar_91}, th. (1.6))  and \cite{MR20:7048}.
\end{demostracion}

\end{document}